\documentclass[leqno,11pt]{amsart}

\newtheorem{Thm}{Theorem}
\newtheorem{Def}{Definition}[section]
\newtheorem{Rem}[Def]{Remark}
\newtheorem{Prop}[Def]{Proposition}
\newtheorem{Cor}[Def]{Corollary}
\newtheorem{Lem}[Def]{Lemma}

\def \R{{\mathbf R}}
\def \T{{\mathbf T}}

\def \wL{w\hbox{-}L}
\def \DD{{\mathcal D}}
\def \V{\hbox{{\rm H}}}
\def \H{\hbox{{\rm V}}}

\def \Ker{\hbox{{\rm Ker}}}

\def \eps{{\varepsilon}}
\def \e{{\varepsilon}}

\def \LL{{\mathcal L}}

\def \OO{{\mathcal O}}
\def \SS{{\mathcal S}}

\def \d{{\partial}}
\def \iint{\int \! \! \! \int }

\def \DOT{\! \cdot \!}

\def \GRAD{\nabla\!}
\def \DIV{\nabla\! \! \cdot }
\def  \ROT{\mbox{rot} \:}
\def \LAP{\Delta}

\newcommand{\WWe}{W_{\eps,3}}
\newcommand{\ue}{ u_\eps}
\newcommand{\ub}{ \overline u_\eps}
\newcommand{\we}{ w_\eps}
\newcommand{\We}{ W_\eps}
\newcommand{\Weu}{ W_{\eps,1}}
\newcommand{\Wed}{ W_{\eps,2}}
\newcommand{\Wet}{ W_{\eps,3}}
\newcommand{\weu}{ w_{\eps,1}}
\newcommand{\wet}{ w_{\eps,3}}
\newcommand{\wedd}{ w_{\eps,2}}
\newcommand {\wde}{w_\eps^\delta}
\newcommand {\wed}{w_\eps^\delta}
\newcommand{\ru}{ \rho_{\eps,1}^\delta}
\newcommand{\rd}{ \rho_{\eps,2}^\delta}
\newcommand{\rt}{ \rho_{\eps,3}^\delta}
\newcommand{\Wd}{ W_\eps^\delta}

\newcommand{\Wdu}{ W_{\eps,1}^\delta}
\newcommand{\Wdd}{ W_{\eps,2}^\delta}
\newcommand{\wdu}{ w_{\eps,1}^\delta}
\newcommand{\wdd}{ w_{\eps,2}^\delta}

\newcommand{\WWd}{ W_{\eps,3}^\delta}
\newcommand{\Wdt}{ W_{\eps,3}^\delta}
\newcommand{\ul}{ \overline u}

\def\eqdefa{\buildrel\hbox{\footnotesize def}\over =}
\newcommand{\Pp}{\Pi_\perp}

\begin{document}
\title[Weak convergence results for rotating fluids]{Weak convergence
   results
\\for inhomogeneous rotating fluid equations }

\author[I. Gallagher]{Isabelle Gallagher}
\address[I. Gallagher]%
{Centre de Math{\'e}matiques UMR 7640 \\
    Ecole Polytechnique \\
     \\
91128 Palaiseau
   \\
    FRANCE}
\email{Isabelle.Gallagher@math.polytechnique.fr}

\author[L. Saint-Raymond]{Laure Saint-Raymond}
\address[L. Saint-Raymond]%
{ Laboratoire J.-L. Lions UMR 7598\\ Universit{\'e} Paris
VI\\
175, rue du Chevaleret\\ 75013 Paris\\FRANCE }
\email{saintray@ann.jussieu.fr }


\maketitle
{\small{\bf Abstract.}
We consider the equations governing incompressible, viscous fluids in
three space dimensions, rotating around an inhomogeneous vector~$ B(x)$: this
is a generalization of the usual rotating fluid model (where~$ B$ is
constant). We  prove the
weak convergence of Leray--type solutions towards a vector field which
satisfies the usual 2D Navier--Stokes equation in the regions of space
where~$ B$ is constant, with Dirichlet boundary conditions, and a
heat--type equation elsewhere. The method of proof uses
weak compactness arguments.}

\vskip1.5cm

\centerline{\bf R\'esultats de convergence faible}
\centerline{\bf
pour des \'equations des fluides tournants non homog\`enes}

\vskip .4cm\noindent 

{\small{\bf R\'esum\'e.}
On consid\`ere les \'equations mod\'elisant des fluides
incompressibles et visqueux  en trois dimensions d'espace, en rotation rapide
autour d'un vecteur non homog\`ene~$ B(x)$: on g\'en\'eralise ainsi le
mod\`ele  habituel des fluides tournants (o\`u~$ B$ est constant).
On montre la convergence des solutions de Leray vers un champ de
vecteurs qui  v\'erifie les \'equations  habituelles de Navier--Stokes
2D dans les r\'egions de l'espace o\`u~$ B$ est constant, avec des
conditions aux  limites de Dirichlet, et une \'equation de type
chaleur  ailleurs. La m\'ethode de d\'emonstration repose sur des
arguments de  compacit\'e faible.}


\section{Introduction}
\setcounter{equation}{0}
The aim of this article is to study the asymptotics of solutions of
rotating fluid equations, in the case when the rotation vector is non
homogeneous. We consider a domain~$ \Omega = \Omega_h\times \Omega_3$,
where~$  \Omega_h$ denotes either the whole space~$\R^2$ or any
periodic domain of~$\R^2$, and similarly~$\Omega_3 $ denotes~$ \R$
or~$ \T$. We are interested in the following system:
\begin{equation}
\label{NS}
\begin{aligned}
\,&\d _t u+u\DOT \GRAD u -\nu \LAP u +\frac1\eps u\wedge
B+\GRAD p=0\quad \hbox{ on } \R^+\times \Omega ,\\
&\DIV u=0\quad \hbox{ on } \R^+\times \Omega,\\
& u_{|t=0}=u^0\quad \hbox{ on }  \Omega
\end{aligned}
\end{equation}
where  $B=b e_3$ is the rotation vector, and~$b$ is a smooth vector
field defined in~$ \Omega_h$. We shall suppose throughout this paper that~$ b$
does not vanish, and is equal to a positive constant~$ b_0$ perturbated by a
smooth~$ C^\infty_c$ function; more assumptions on~$ b$ will be made as we go
along. Before stating the result we shall prove here, let us
recall some well-known facts in the constant case ($ b = 1$). The
rotating fluid equations, with~$ b$ constant and homogeneous,
modelize the movement of
the atmosphere or the oceans at mid-latitudes (see for
instance~\cite{Greenspan} or~\cite{Pedlosky}). The fluid is supposed
to be incompressible (which corresponds to the hydrostatic
approximation), and its viscosity is~$ \nu > 0$. The vector field~$ u$
is the velocity and the scalar~$ p$ is the pressure, both are
unknown.
The parameter~$ \eps$ is the Rossby number, and its inverse stands for
the speed of rotation of the Earth. Taking the limit~$ \eps \to 0$
means that the scale of motion of the fluid is much smaller than that
of the Earth.
Note that one can also see~$ B$ as a magnetic field, in which case it makes
sense to understand what happens when~$b $ is not homogeneous; that also
holds if one  wants to study the movement of the atmosphere in
other regions than mid-latitudes.

In the constant case, those equations have been studied by a number of
authors. We
refer for instance to the works of A. Babin, A. Mahalov and
B. Nicolaenko~\cite{basil}-\cite{basil3}, I. Gallagher~\cite{I3},
E. Grenier~\cite{Grenier} for the periodic case, and J.-Y. Chemin,
B. Desjardins, I. Gallagher and E. Grenier~\cite{cdgg} for the whole
space case as well as~\cite{cdggplack} for the case of horizontal
plates with Dirichlet boundary
conditions (for such boundary conditions we refer also
to the work of E. Grenier and N. Masmoudi~\cite{greniermasmoudi} as
well as N. Masmoudi~\cite{masmoudi}). The results in those papers 
concern both weak and strong
solutions; in
this article we shall only be concerned with Leray--type weak
solutions (\cite{leray}): we will see in Section~\ref{sing} below that
their existence is an easy adaptation of the proof of Leray's
existence theorem~\cite{leray}. In
the constant case, it is known that weak solutions converge towards
the solution of the
two--dimensional Navier-Stokes equations. Such a result in the whole
space case is due to Strichartz-type estimates (which are obtained by
writing  the solution of the linearized problem in Fourier space),
whereas in the
periodic case it follows from the study of the (discrete) spectrum of
the rotating fluid operator (following methods introduced by
S. Schochet in~\cite{schochet}).

The problem here if we want to follow those methods is that it
does not seem a good idea to take the Fourier transform of the operator
$$
L u \eqdefa P(u \wedge B), \quad \nabla \cdot u = 0,
$$
where $P$ denotes the Leray projection onto divergence free vector fields,
when~$ B$ is not homogeneous; moreover the study of the spectrum of~$
L$ is not an easy matter. So our strategy to study this
problem is first to try and recover the well-known results of the
constant case without using any information on the spectrum of~$ L$
(other than the determination of its kernel),
and without using the Fourier transform. This will be achieved in
Section~\ref{constant}. Then the study of the variable case will be an
adaptation of the constant case, in Section~\ref{variable}.

Before stating the results we shall prove in this paper, let us
comment on the difficulties compared with the constant case: as stated
above, it is easy to construct a bounded family of weak solutions to
our problem, whether~$ b$ is constant or not. Hence one can construct
a weak limit point~$ \overline u$, and the question we want to address
is to find the equation satisfied by~$ \overline u$. Of course the
problem consists in taking the limit in the non linear part of the
equation. As noted above, we do not wish to study the spectrum of the
operator~$L $ since that seems to be a difficult issue. So we cannot
apply the usual, constant~$ b$ methods, as to our knowledge they all
involve spectral properties of~$ L$. The idea therefore is to turn to
what is known as ``weak compactness methods'', in the spirit of Lions
and Masmoudi~\cite{lionsmasmoudi1}-\cite{lionsmasmoudi2} (for the 
incompressible limit). We
shall recall briefly below what those methods are, and then we shall
state the main results of this paper.

\subsection{ Weak compactness methods}

Let us explain what weak compactness methods are all
about. The idea is as follows: as usual the trouble to find the limit
of the equation comes from the bilinear terms. They can be  separated
into three categories:

\begin{itemize}
\item products involving only elements of the kernel of the
penalization~$ L$, which can be
shown to be compact (see Corollary~\ref{t-reg});

\item  products of elements of the kernel against
    elements of~$(\Ker L)^\perp$, for which one can take the limit since
    elements of~$(\Ker L)^\perp$ converge weakly to zero (see
Corollary~\ref{t-reg});

\item products  involving only elements of~$(\Ker L)^\perp$, which are
    the problem.
\end{itemize}

The idea now is to prove that in the last situation, the limit is in
fact zero for
   algebraic reasons: in previous works on rotating fluids,
that result was proved essentially by writing the product of two elements
of~$(\Ker L)^\perp$ by projection onto eigenvectors of~$ L$. In the
   periodic case, a ``miracle'' in the formulation yielded the result
   (see~\cite{basil}-\cite{basil3} or~\cite{I3}),
   whereas in the whole space case, Strichartz estimates did the job
   (and the convergence was strong), see~\cite{cdgg}. In this paper we
   will show that the result has in fact not much to do with spectral
   properties of~$ L$, but is due to simple algebraic properties. Let 
us recall the result in the case of the
   incompressible limit, where such properties were first used
(see~\cite{lionsmasmoudi1}).

\begin{Prop}\cite{lionsmasmoudi1}
{\sl
Let $(\rho_\eps),(u_\eps),(\theta_\eps)$ be
bounded families of $L^2([0,T],H^1(\Omega))$ such that
$$\rho_\eps \rightharpoonup \rho,\quad u_\eps \rightharpoonup u,\quad
\theta_\eps
\rightharpoonup \theta\hbox{ as }\eps \to 0.$$
Assume that
$$\begin{array}l
\displaystyle \eps\d_t \rho_\eps +\DIV u_\eps =0,\\
\displaystyle \eps\d_t u_\eps +\GRAD(\rho_\eps +\theta_\eps)=s_\eps,\\
\displaystyle\eps \d _t \theta_\eps+\frac23 \DIV u_\eps=s'_\eps,
\end{array}
$$
where $ s_\eps,s'_\eps\to 0$ in $L^1([0,T],H^{-s}(\Omega))$ for some
$s>0$. Then
$$P\DIV(u_\eps\otimes u_\eps)\to P\DIV(u\otimes u) \hbox{ and } \DIV(u_\eps
\theta_\eps) \to \DIV(u\theta)$$
in the sense of distributions.
}
\end{Prop}
\begin{proof}
This result has to be compared with the so-called ``compensated compactness"
theorems, in the sense that the convergences of some quadratic quantities in
$\rho_\eps,u_\eps,\theta_\eps$ are established under the assumption that some
combinations of the derivatives of these functions converge strongly in time to
0.
The proof consists in checking that the acoustic oscillations  do not
bring any contribution to the limiting terms. We introduce the
following decompositions:
$$
\ue=P\ue+\nabla \psi_\eps,\quad \hbox{ and }\quad\theta_\eps=
{3\theta_\eps-2\rho_\eps\over 5} +\pi_\eps,
$$
so that
$$\begin{aligned}
  P\ue \: \hbox{ and } \: 
{3\theta_\eps-2\rho_\eps\over 5}  \: \hbox{ are bounded 
in
}W^{1,1}([0,T],H^{-s}(\Omega)), \quad \mbox{and}\\
\eps \d _t 
\nabla \psi_\eps +\nabla \pi_\eps = (Id - P)s_\eps\to 0 
,\quad
\eps
\d _t
\pi_\eps+{2
\over 3}\Delta \psi_\eps= \frac{2}{5} 
s'_\eps \to 0\hbox{ in 
}
L^1([0,T],H^{-s}(\Omega)).
\end{aligned}$$
We shall note in the 
following~$ S_\eps \eqdefa (Id - P)s_\eps$
and~$\displaystyle 
S'_\eps \eqdefa \frac{2}{5} s'_\eps$.  
The incompressibility and 
Boussinesq relations
$$\DIV u=0,\quad \GRAD(\rho+\theta)=0$$
allow to 
identify the limits
$$\begin{aligned}
P\ue\to u,\quad 
{3\theta_\eps-2\rho_\eps\over 5}\to \theta &\hbox{ 
in
}L^2([0,T]\times \Omega),\\
\nabla \psi_\eps\rightharpoonup 
0,\quad \pi_\eps \rightharpoonup 0 &\hbox{
in }w-L^2([0,T]\times 
\Omega),
\end{aligned}$$
from which we deduce that, in the sense of 
distributions
$$\begin{aligned}
\ue\otimes \ue-\nabla 
\psi_\eps\otimes \nabla \psi_\eps&\to u\otimes 
u,
\\
\theta_\eps\ue-\pi_\eps \nabla \psi_\eps &\to \theta 
u.
\end{aligned}$$
The key argument is therefore the following formal 
computation (which 
can be made
rigorous by introducing 
regularizations with respect to the space variable 
$x$)
$$\begin{aligned}
P\nabla.(\nabla \psi_\eps \otimes
\nabla 
\psi_\eps)
&={1\over 2}P\nabla | \nabla 
\psi_\eps|^2+P(\Delta
\psi_\eps \nabla \psi_\eps) \\
&= {3\over 2} 
P(-\d _t (\eps \pi_\eps \nabla
\psi_\eps)-\pi_\eps \nabla \pi_\eps 
+
\pi_\eps S_\eps +{S'}_\eps \nabla \psi_\eps)
\cr
&={3\over 2} P(-\d 
_t (\eps \pi_\eps \nabla
\psi_\eps) +
\pi_\eps S_\eps +{S'}_\eps 
\nabla \psi_\eps),
\\
\nabla.( \pi_\eps\nabla
\psi_\eps)
&= 
\pi_\eps
\Delta \psi_\eps+\nabla \psi_\eps.\nabla
\pi_\eps\\
&= 
{3\over 2}\pi_\eps ({S'}_\eps-\eps\d 
_t
\pi_\eps)+\nabla
\psi_\eps.(S_\eps-\eps \d _t \nabla \psi_\eps
)\\ 
&={3\over 2}\pi_\eps {S'}_\eps+\nabla
\psi_\eps.S_\eps-{3\eps\over 
4}\d _t
|\pi_\eps|^2-{\eps\over 2} \d _t |\nabla
\psi_\eps|^2,
\end{aligned}$$
which shows that the contribution of the acoustic 
oscillations is
negligible.
\end{proof}
Inspired by the previous 
computation, we shall in this article try to
use a similar method in 
the case of rotating fluids: we refer to the
proofs 
of
Propositions~\ref{couplage} and~\ref{couplagevariable} for precise 
computations.

\subsection{ Main results}

Since we consider 
incompressible flows,  we introduce
  the following subspaces of 
$L^2(\Omega)$ and $H^1(\Omega)$
$$\V=\{u\in L^2(\Omega)\,/\, \DIV 
u=0\},\quad \H=\{u\in H^1(\Omega)\,/\, \DIV
u=0\}.$$
We will also use 
the following notation for the inhomogeneous Sobolev spaces
$$ 
H^s(\Omega)=\{u \in \DD'(\Omega)\,/\, (Id -\LAP)^{s/2} u\in 
L^2(\Omega) \}.$$
Similarly homogeneous Sobolev spaces will be 
defined by
$$\dot H^s(\Omega)=\{u \in \DD'(\Omega)\,/\, (-\LAP)^{s/2} 
u\in 
L^2(\Omega) \}.$$
  It will appear clearly in the following 
that the horizontal variables
play a special role in this problem. 
Consequently we shall use the following
notation: if~$ x$ is a point 
in~$
\Omega$, then we shall note its cartesian coordinates 
by~$(x_1,x_2,x_3) $, and
the horizontal part of~$ x$ will be 
denoted~$ x_h \eqdefa (x_1,x_2)
\in \Omega_h$. Similarly
we will 
denote the horizontal part of any vector  field~$f$ by~$f_h$, 
the
horizontal gradient by~$
\nabla_h
\eqdefa (\d_1,\d_2)$ and its 
orthogonal by~$\nabla_h^\perp = (\d_2,-\d_1)$, and
the horizontal 
divergence and Laplacian respectively by~$ \mbox{div}_h
f \eqdefa 
\d_1 f_1 + \d_2 f_2 = \nabla_h \cdot f_h$ and~$ \Delta_h
\eqdefa 
\d_1^2 + \d_2^2$.

Finally as usual, $ C$ will denote a constant 
which can change from
line to line, and~$\nabla p $ will denote the 
gradient of a function
which can also change from
line to line.

Before stating the main theorems of this paper, let us give some
additional definitions. We will note by
$$
S\eqdefa \{x\in \Omega\,/\, \GRAD b(x)=0\} \quad \mbox{and} \quad \OO
\eqdefa  \{x\in \Omega\,/\, \GRAD b(x)\neq0\}.
$$\
Finally~$ \SS$ will be the interior of the singular set~$ S$. We will
assume in the sequel that
$$
\hbox{ The set } \Omega \setminus \left(\SS \cup \OO \right)   \hbox{
   is of Lebesgue measure 0 ;}\leqno(H0)
$$
$$
\SS \hbox{ is a
smooth domain,}\leqno(H1)
$$
$$
\begin{array}{c}
\hbox{ On each connected component }  \OO_j \hbox{  of } \OO,  \hbox{ 
there is a smooth
function } \sigma_j \\
\hbox{ such that } (b,\sigma_j,x_3) \hbox{ is a global
   smooth coordinate system  } \\
\hbox{ and }\OO_j = \{ x_h \in \R^2 \: | \:
\left (b(x_h),\sigma_j (x_h)\right) \in B_j \times \Sigma_j )\} 
\times \Omega_3.
\end{array} \leqno(H2)
$$

Now we are ready to state the main theorems of this paper. The first
result, rather standard, shows that there are weak solutions to the 
system~(\ref{NS}).

\begin{Thm}\label{thm:energie}
{\sl
Let $u^0$ be any  vector field in $\V$.
Then for all $\eps>0$, Equation (\ref{NS}) has at least
one weak solution $\ue\in L^\infty(\R^+,\V)\cap  L^2(\R^+,\dot H^1)$.
Moreover, for all $t>0$, the following energy estimate holds:
\begin{equation}
\label{eq:energie}
\|\ue(t)\|^2_{L^2} +2\nu \int_0^t \|\GRAD \ue(s)\|_{L^2}^2
ds\leq
\|u^0\|_{L^2}^2.
\end{equation}
}
\end{Thm}

Now the aim of the paper is to describe the limit of~$\ue$ as~$ \eps$
goes to zero.
We will first concentrate on the constant case.

\begin{Thm}\label{thm:constant}
{\sl
Suppose that~$ B = b e_3$ where~$ b$ is constant and homogeneous.
Let $u^0$ be any  vector field in $\V$, and let~$ \ue$ be any weak
solution of~(\ref{NS}) in the sense of Theorem~\ref{thm:energie}. 
Then~$\ue $ converges weakly
in~$ L^2_{loc}(\R^+\times \Omega)$
  to a limit~$ \overline u$ which if~$ \Omega_3 = \R$ is zero, and if~$
  \Omega_3 = \T$ is the solution of the two dimensional Navier--Stokes
  equations
$$
(NS2D) \quad \quad \partial_t \overline u - \nu \Delta_h \overline u 
+  \overline u_h \cdot
\nabla_h \overline u = (-\nabla_h p, 0), \quad \mbox{div}_h  \overline u_h =
0, \quad  \overline u_{|t = 0} =
\int_{\T} u^0(x_h, x_3) \:dx_3.
$$
}
\end{Thm}
\begin{Rem}
{\sl
This theorem is by no means a novelty, it is even rather less precise
than other such results one can find in the literature
(\cite{basil}-\cite{basil3}, \cite{cdgg},  \cite{I3},
\cite{Grenier}). As  we will see in Section~\ref{constant}, the 
interest of this result lies in its
proof, which contrary to the references above, does not depend on the
boundary conditions (which can be the whole space or periodic, in each
direction).
}
\end{Rem}

Now let us state the new result of this paper, concerning the case
when~$ b$ is not homogeneous.

\begin{Thm}\label{thm:variable}
{\sl
Suppose that~$ B = b e_3$ where~$ b = b(x_h)$ is a nonnegative smooth function,
say a~$C^\infty_c (\Omega_h)$ perturbation of a constant, and where 
assumptions~$(H0) $
to~$(H2) $ are satisfied.
Let $u^0$ be any  vector field in $\V$, and let~$ \ue$ be any weak
solution of~(\ref{NS})  in the sense of Theorem~\ref{thm:energie}. 
Then~$\ue $ converges weakly
in~$ L^2_{loc}(\R^+\times\Omega)$
  to a limit~$ \overline u$ which if~$ \Omega_3 = \R$ is zero, and if~$
  \Omega_3 = \T$ is defined as follows: the third component~$\overline 
u^3 $ satisfies the transport equation
$$
\partial_t \overline u^3  - \nu \Delta_h \overline u^3 +  \overline u_h \cdot
\nabla_h \overline u^3 = 0, \quad \partial_3 \overline u^3 = 0, \quad 
\overline u^3_{|t = 0} =
\int_{\T} u^{0,3}(x_h, x_3) \:dx_3 \quad \mbox{in} \quad \R^+ \times 
(\OO \cup \SS),
$$
while the horizontal component~$\overline u_h $ depends on the region
of space considered:
\begin{itemize}

\item in $ \SS$, $ \overline u_h$
  satisfies the two dimensional  Navier--Stokes
  equations~$ (NS2D)$ with Dirichlet boundary conditions.

\item in  $ \OO$,  $ \overline u_h $ satisfies the following heat equation:
$$
\partial_t \overline u_h - \nu \Pi \Delta_h \overline u_h =0,
$$
where~$ \Pi$ is the~$ L^2$ orthogonal projection onto the kernel of~$
L$ (which can be extended to~$ {\mathcal D}'(\OO \cup \SS)$) which 
satisfies in particular
$$
- \left( \Pi \Delta_h \overline u_h |  \overline u_h\right)_{L^2(\OO)} =
\|\nabla_h \overline u_h\|_{L^2(\OO)}^2.
$$
\end{itemize}
}
\end{Thm}
\begin{Rem}
{\sl 
 In the regions 
where~$ b$
is homogeneous, we recover at the limit the 2D 
Navier--Stokes
equations as usual. The Dirichlet boundary conditions 
appear quite
naturally, considering that on the other side of the 
boundary one
finds that~$ \overline u_h$ is proportional 
to~$\nabla_h^\perp b $
which vanishes on the boundary of~$ \SS$. More 
surprising is certainly
the result in the region where~$ b$ is not 
homogeneous.  This
can be understood as some sort of turbulent 
behaviour, where all
scales are mixed due to the variation of~$ b$. 
Technically the result
is due to the fact that the kernel of~$ L$ is 
very small as soon as~$
b$ is not a constant, which induces a lot of 
rigidity in the limit equation.
}
\end{Rem}

The structure of the paper is as follows. In the next 
section, we
present the operator~$ L$ and study its main properties 
(proof of
Theorem~\ref{thm:energie}, kernel,
projections onto~$ \Ker 
L$). The following section is devoted to the
proof of 
Theorem~\ref{thm:constant}. Although the result is not new, 
we
present an alternative proof which holds regardless of the 
domain
(with no boundary). This serves as a warm--up to the final 
section, in
which the general variable case is presented, with the 
proof of Theorem~\ref{thm:variable}.

\begin{Rem}
{\sl 
 One can 
wonder about what remains of those results under more general 
assumptions on the rotation vector~$ B$. First let us consider the
case when~$ B = b(t,x_h) e_3$
 depends also on time. The results  of 
Section~\ref{sing} are
 identical, and if the sets~$
 \SS$ and~$ \OO$ 
are independent of time one recovers the same type
 of theorem as in 
the constant case. In particular  the equation on~$ \overline
 u$ 
in~$ \OO$ is derived in an identical way to Section~\ref{oo}. If the 
sets~$
 \SS$ and~$ \OO$
 do depend on time, then one  has to be a 
little bit more
 careful and this issue will not be treated here.
A 
more physical problem is the case when the direction of~$ B$ is 
not
fixed, in other words when~$ B$ is a three component vector, 
depending on
all three variables.   Then geometrical problems appear, 
simply to
determine the kernel of~$ L$~; this will be dealt with in a 
forecoming paper.
}
\end{Rem}


\section{Study of the singular perturbation}
\label{sing}
\setcounter{equation}{0}
\subsection{Energy estimate}
In this section we shall prove Theorem~\ref{thm:energie} stated in 
the introduction.
\begin{proof}
The structure of the equation (\ref{NS}) governing the
rotating fluids is very similar to the one of the usual
Navier-Stokes equation, since the singular perturbation is
just a linear skew-symmetric operator. Therefore weak
solutions `` {\`a} la Leray " can be constructed by the
approximation scheme of Friedrichs~:
approximate solutions are obtained by a standard truncation~$ J_n$
of high frequencies. In order to obtain uniform bounds on
these approximate solutions, we have just to check that the
energy inequality is still satisfied.
Computing formally the $L^2$ scalar product of (\ref{NS}) by
$u$ leads to
$$\frac12 {d\over dt} \|u\|_{L^2}^2 =-\int\!\! \left(\frac12
(u\DOT
\GRAD) |u|^2-\nu u\DOT\LAP u +\frac1\eps u\DOT u
\wedge B+ u\DOT \GRAD p\right) \!dx.$$
Integrating by parts (without boundary) and using the
incompressibility constraint, we get
$$\frac12 {d\over dt} \|u\|_{L^2}^2 =-\nu \|\GRAD
u\|_{L^2}^2,$$
which holds for any smooth solution of (\ref{NS}).

The energy inequality for weak solutions is obtained by
taking limits in the approximation scheme.
\end{proof}

\smallskip
In particular, the  energy estimate provides uniform
bounds in
$L^\infty(\R^+,\V) \cap L^2(\R^+,\dot H^1)$
on any family $(\ue)_{\eps>0}$ of weak solutions of (\ref{NS})
provided that the initial data $u^0$ belongs to $\V$.

\begin{Cor}\label{cv-faible}
{\sl
Let $u^0$ be any vector field in $\V$.
For all $\eps>0$, denote by $\ue$ a weak solution of
(\ref{NS}). Then there exists $\ul \in
L^\infty(\R^+,\V) \cap L^2(\R^+,\dot H^1)$, such
that, up to extraction of a subsequence,
$$\ue \rightharpoonup \ul \hbox{ in }\wL^2_{loc}(\R^+\times
\Omega) \hbox{ as }
\eps
\to 0.$$
}
\end{Cor}

\subsection{ Characterization of the kernel}

We are 
interested in describing the asymptotic behaviour of
$(\ue)$, i.e. in 
characterizing its limit points. Of course,
the equations satisfied 
by such a limit point $\ul$ depend
strongly on the structure of the 
singular perturbation
\begin{equation}
\label{L-def}
L :u\in \V 
\mapsto P(u\wedge B)\in
\V
\end{equation}
  where $P$ denotes the 
Leray projection from
$L^2(\Omega)$ onto its subspace $\V$ of
divergence-free vector fields. In particular, we will prove
that $\ul 
$ belongs to the kernel $\Ker(L)$ of $L$,
which is characterized in 
the following proposition.

\begin{Prop}\label{noyau}
{\sl
Define the 
linear operator $L$ by (\ref{L-def}).  Then
$u\in \V$ 
belongs to $\Ker(L)$ if and only
if there exist  $\GRAD_h \varphi 
\in
 L^2(\Omega_h)$ and $\alpha \in
L^2(\Omega_h)$ with
$$\GRAD_h b\DOT \GRAD^\perp_h \varphi=0, $$
such that
$$u= \GRAD^\perp_h \varphi+\alpha e_3.$$
}
\end{Prop}

\begin{proof}
We have
$$P(u\wedge B)=0.$$
Then, in the sense of distributions,
$$\ROT(u\wedge B)=0,$$
which can be rewritten
$$(\DIV B)u+(B\DOT \GRAD)u-(u\DOT \GRAD)B-(\DIV u) B=0.$$
As $\DIV B=\DIV u=0$ and $B=b e_3$, we get
\begin{equation}
\label{kernel}
b\d _3 u-(u\DOT \GRAD)b e_3=0.
\end{equation}

In particular, $\d_3 u_1=\d_3 u_2=0$ from which we deduce that
\begin{equation}
\label{x3-dep}
u_1,u_2\in L^2(\Omega_h).
\end{equation}
Note that in the case where $\Omega_3=\R$, the invariance with
respect to~$x_3$ and the fact that~$u\in L^2(\Omega)$ imply
that~$u_1=u_2=0$ (and therefore~$ u_3 = 0$ by the divergence free condition).

Differentiating the incompressibility constraint with respect to $x_3$
leads then to
$$\d^2_{33} u_3=-\d^2_{13}u_1-\d^2_{23}u_2=0$$
in the sense of distributions. The function $\d_3 u_3$ depends only on $x_1$
and $x_2$, and satisfies~$\displaystyle \int \d_3 u_3 dx_3=0$. 
So~$\d_3 u_3=0$ and
\begin{equation}
\label{divh}
u_3\in
L^2(\Omega_h),\quad \d_1 u_1+\d_2 u_2=0.
\end{equation}
Combining (\ref{x3-dep}) and (\ref{divh}) provides the existence of
$\GRAD _h\varphi \in
L^2 (\Omega_h)$ such that
$$u_1=\d_2\varphi,\quad u_2=-\d_1 \varphi.$$

Replacing in (\ref{kernel}) leads to
$$\GRAD_h^\perp \varphi\DOT \GRAD_h b=b\d_3 u_3 =0,$$
which concludes the proof.
\end{proof}

\smallskip
Before applying this result to the characterization of the weak limit $\ul$,
let us just specify it in two important cases. If $\GRAD
b=0$ almost everywhere, $u\in
\V$ belongs to $\Ker(L)$ if and only
if
$$
u= \GRAD^\perp _h\varphi+\alpha e_3,
$$
for some $\GRAD_h\varphi \in
L^2(\Omega_h)$ and $\alpha \in
L^2(\Omega_h)$.
   If $\GRAD
b\neq 0$ almost everywhere (in other words, if~$ \Omega \setminus \OO$
   is of Lebesgue measure zero), then the condition arising on $u$ is much more
restrictive~:
$u\in \V$ belongs to $\Ker(L)$ if and only
if on each connected component of~$\OO $,
$$u= F(b) \GRAD^\perp_h b+\alpha e_3,$$
for some square integrable function~$ F(b)$ and some~$\alpha \in 
L^2(\Omega_h)$.

\bigskip
  From this characterization of $\Ker(L)$, we deduce
some constraints on the weak limit $\ul$.

\begin{Cor}\label{contrainte}
{\sl
Let $u^0$ be any vector field in $\V$.
Denote by $(\ue)_{\eps>0}$ a family of weak solutions of
(\ref{NS}), and by $\ul$  any
of its limit points. Then, there exist $\varphi \in L^\infty(\R^+,\dot
H^1(\Omega_h)) \cap L^2(\R^+,\dot H^2(\Omega_h))$ and $\alpha \in
L^\infty(\R^+,L^2(\Omega_h)) \cap L^2(\R^+,\dot H^1(\Omega_h))$ with
$$\GRAD_h b\DOT \GRAD_h^\perp \varphi=0, $$
such that
$$\ul= \GRAD_h^\perp \varphi+\alpha e_3.$$
}
\end{Cor}

\begin{proof} Let $\chi\in \DD(\R^+ \times\Omega) $ be any
divergence-free test function. Multiplying (\ref{NS}) by
$\eps \chi$ and integrating with respect to all variables
leads to
$$ \iint \ue(\eps \d _t \chi+\eps \ue\DOT \GRAD  \chi+\eps \nu
\LAP \chi +\chi\wedge B)dxdt=0.$$
Because of the bounds coming from the energy estimate, we can
take limits in the previous identity as $\eps \to 0$ to get
$$\int \ul \wedge B\DOT \chi dxdt=0.$$
This means that there exists some $p$ such that
$$\ul \wedge B=\GRAD p.$$
As $\ue$ satisfies the incompressibility relation for all
$\eps>0$,
$$\DIV
\ul =0.$$
Then $\ul(t) \in \Ker(L)$  for almost all $t\in \R^+$, and we conclude by
Proposition \ref{noyau} that
   there exist $\varphi \in L^\infty(\R^+,
\dot H^1(\Omega_h))\cap L^2(\R^+,
\dot H^2(\Omega_h))$ and $\alpha \in L^\infty(\R^+, 
L^2(\Omega_h))\cap L^2(\R^+,
\dot H^1(\Omega_h))$ with
$$\GRAD_h b\DOT \GRAD_h^\perp \varphi=0, $$
such that
$$u= \GRAD_h^\perp \varphi+\alpha e_3.$$
\end{proof}

\subsection{ Decomposition by projection on the kernel}
   To go further in the description of the asymptotic behaviour
(i.e. in the characterization of $\ul$), we have to isolate
the fast oscillations generated by the singular perturbation
$L$, which produce ``big" terms in (\ref{NS}), but converge
weakly to 0.

Therefore we introduce the following decomposition
$$\ue=\ub +\we,$$
where $\ub=\Pi \ue$ is the $ L^2$ orthogonal projection of $\ue$ on $\Ker(L)$
and $\we=\Pp \ue$ is the projection of $\ue$ on~$\Ker(L)^\perp$.

We have seen in the previous paragraph that the characterization of the kernel
$\Ker(L)$ is strongly linked to the geometry of the vector field
$b$.
In
order to obtain further regularity properties on $\Pi$ and $\Pp$, we then
need a precise description of the singular set
$$S= \{x\in \Omega\,/\, \GRAD b(x)=0\},$$
which justifies  Assumptions~$ (H0)$ to~$ (H2)$ given in the introduction.

Before
stating the main properties of
$\Pi$ and $\Pp$, let us give the following definitions: by
Assumption~$(H0)$ it is enough to describe the limiting function~$
\overline u$ on~$ \OO \cup \SS$. So it is natural to define the
following  function spaces: for all~$ s \geq 0$, $  H^s (\OO \cup \SS)
$ is the closure of~$ C^\infty_c (\OO \cup \SS)$ for the~$ H^s$ norm,
and we will note, for~$ s \geq 0$, $  H^{-s} (\OO \cup \SS)$   the
dual space of~ $  H^s (\OO \cup \SS)$.

It will be useful in the following to note that
\begin{equation}
\label{eq:inclusion}
\forall s \geq 0, \quad H^s (\OO \cup \SS) \subset  H^s (\Omega)
\quad \mbox{and} \quad \forall s \leq 0, \quad H^s (\Omega) \subset 
H^s (\OO \cup \SS).
\end{equation}
\begin{Prop}\label{decomposition}
{\sl
Define the linear operator $L$ by (\ref{L-def}). Denote by
$\Pi$ the orthogonal projection of $\V$ onto
$\Ker(L)$ and by $\Pp=Id-\Pi$ the orthogonal projection of $\V$ onto
$\Ker(L)^\perp$. The operators $\Pi$ and $\Pp$ so
defined
have
natural extensions to tempered distributions on $\OO \cup \SS$,
and for all
$s\in \R$, there exists some
$C_s>0$ such that for any function~$ u\in H^s(\OO \cup \SS) $,
$$
\| \Pi u\|_{H^s(\OO \cup \SS)} \leq
C_s\|u
\|_{H^s(\OO \cup \SS)}
\quad
\mbox{and} \quad \| \Pp u\|_{H^s(\OO \cup \SS)} \leq
C_s\|u
\|_{H^s(\OO \cup \SS)}.
$$
}
\end{Prop}
  \begin{proof}
By Proposition \ref{noyau}, for all $u\in \V$
$$\Pi u=\GRAD^\perp_h \varphi +\alpha e_3,$$
for some $\GRAD_h\varphi \in
   L^2(\Omega_h)$ and $\alpha \in L^2(\Omega_h)$ with
$$\GRAD_h b\DOT \GRAD_h^\perp \varphi=0. $$
By definition $\Pp u=u-\Pi u$ is orthogonal to any element of $\Ker(L)$.
So for all $\beta \in L^2(\Omega_h)$,
$$\int (u-\Pi u) \DOT \beta e_3 dx=0,$$
which implies that
\begin{equation}
\label{alpha-def}
\alpha = {1\over |\Omega_3|}\int u_3 dx_3.
\end{equation}
In order to determine $\varphi$, we consider separately the domains $\OO$ and
$S$.

On $S$, $\varphi$ is  defined as
$$
\displaystyle
  \GRAD_h^\perp \varphi =\left({1\over |\Omega_3|}\int_{\Omega_3} u_h
dx_3\right),
$$
and the smoothness properties stated in 
Proposition~\ref{decomposition} are obvious.

On~$ \OO$, we use Assumption~$ (H2)$ which  implies that on each
connected component~$ \OO_j$ of~$ \OO$, $ u$ can be written~$u(x) =
f_j(b,\sigma_j,x_3)$ where~$ f_j$ has the same smoothness as~$ u$ since the
change of coordinate is in~$ C^\infty (\Omega)$. Then clearly the
projection~$ \Pi$ is simply defined by
$$
\Pi u \eqdefa \frac{1}{|\Omega_3| |\Sigma_j| }\int_{\Sigma_j}\int_{\Omega_3}
f_j(b,\sigma_j,x_3 ) d\sigma_j dx_3 ,
$$
and the result follows.
\end{proof}

\begin{Cor}\label{t-reg}
{\sl
Let $u^0$ be any vector field in $\V$.
Denote by $(\ue)_{\eps>0}$ a family of weak solutions of
(\ref{NS}), and by $\ul$  any
of its limit points. Consider a subsequence of $(\ue)$
(abusively denoted~$(\ue)$) such that
$$\ue \rightharpoonup \ul \hbox{ in }\wL^2_{loc}(\R^+\times
\Omega) \hbox{ as }
\eps
\to 0.$$
Then, if $\ub=\Pi \ue$ is the projection of $\ue$ on
$\Ker(L)$,
$$\ub \to \ul \hbox{ strongly in
}L^2_{loc}(\R^+\times
(\OO \cup \SS)) \hbox{ as }
\eps
\to 0.$$
}
\end{Cor}

\begin{proof}
By Proposition \ref{decomposition}, the projection $\Pi$
is continuous in~$ L^2$. Then, by the energy estimate,
$$\ub \hbox{ is uniformly bounded in
} L^2_{loc}(\R^+, {\mathcal K})$$
for some compact subset~$ {\mathcal K}$ of~$ L^2(\Omega)$,
which provides regularity with respect to space variables.

The second step consists in getting regularity with respect
to time. Apply $\Pi$ to the convection equation in
(\ref{NS})~:
\begin{equation}
\label{moy}
\d _t \ub+\Pi ( \ue\DOT \GRAD \ue-\nu \LAP \ue)=0.
\end{equation}
As $\ue$ is divergence-free,
$$\ue\DOT \GRAD \ue=\DIV (\ue \otimes \ue)\hbox{ is uniformly
bounded in } L^\infty(\R^+,W^{-1,1}(\Omega)).$$
  From the energy estimate, we also deduce that
$$ \LAP \ue \hbox{ is uniformly bounded in
} L^2(\R^+, H^{-1}(\Omega)),$$
from which we infer that
$$
\ue\DOT \GRAD \ue-\nu \LAP \ue  \hbox{ is
uniformly bounded in }  L^2_{loc}(\R^+,H^{-5/2}_{loc}(\Omega)).
$$
Combining this with~(\ref{eq:inclusion}),  we get by Proposition
\ref{decomposition}
$$\d _t \ub=-\Pi( \ue\DOT \GRAD \ue-\nu \LAP \ue) \hbox{ is
uniformly bounded in } L^2_{loc}(\R^+,H^{-5/2}_{loc}(\OO \cup \SS)),$$
which provides the expected regularity in $t$.

Aubin's lemma~\cite{aubin} then gives the following interpolation result
$$(\ub) \hbox{ is strongly compact in } L^2_{loc}(\R^+\times
(\OO \cup \SS)).$$
By Proposition \ref{decomposition}, $\Pi$ belongs to
$C(L^2(\Omega),L^2(\Omega))$ and therefore to  $C(w-L^2(\Omega),
w-L^2(\Omega))$, from which we deduce that
$$\ub=\Pi \ue \rightharpoonup \Pi \ul=\ul.$$
Combining both results shows that
$\ub$ converges strongly to $\ul$ in $L^2_{loc}(\R^+\times
(\OO \cup \SS)).$
\end{proof}

\subsection{ Remarks concerning the regularity}
\subsubsection{Comparison with the gyrokinetic approximation}
\label{gyro}
As mentionned in the introduction, the study of the asymptotics
for an  inhomogeneous penalization is a natural question in the
magnetohydrodynamic framework, when $B$ represents the magnetic 
field. Such a study has
been performed for the gyrokinetic approximation \cite{GSR}, that is 
for a kinetic model
perturbed by a singular magnetic constraint~:
\begin{itemize}
\item
in the case where $B = b(x_h) e_3$, the singular limit is exactly the 
same as in
the constant case~: the fast rotation has an averaging effect in the 
plane orthogonal to
the magnetic lines;
\item
in the case where $B$ has constant modulus but variable direction, 
extra drift terms are
obtained due to the curvature of the field.
\end{itemize}
A simplified version of this result can be written as follows.
\begin{Thm}\cite{GSR}
{\sl
Let $f^{0}$ be a function of $ L^\infty(\Omega\times \R^3)$, and 
$(f_\eps)$ be a
family of  solutions of
$$\d _t f_\eps +v\cdot \GRAD_x f_\eps +\frac1\eps v\wedge B\cdot 
\GRAD _v f_\eps =0,
\quad t\in
\R^*_+,\, (x,v)\in \Omega\times \R^3,$$
with initial condition
$$ f_\eps(|t=0 )=f^{0}.$$
Then the family $(f_\eps)$ is relatively compact in 
$w*-L^\infty(\R_+\times \Omega\times
\R^3)$, as well as the family $(g_\eps)$ defined by
$$g_\eps (t,x,w)=f_\eps(t,x,R(x,-\frac t\eps)w)$$
where $R(x,\theta)$ denotes the rotation of  angle $\theta$ around 
the oriented axis of
direction $B(x)$. Moreover,
\begin{itemize}
\item if $B=b e_3$ with $b\in C^1(\Omega_h,\R_+^*)$, any  limit point 
of $(g_\eps)$
satisfies
$$ \d _t g+v_3 \d _{x_3} g=0;$$
\item if $B\in C^1(\Omega)$ with $\nabla_x \cdot B=0$ and $|B|\equiv 
1$, any of its limit
points satisfies
$$\d_t g+ (w\cdot B)B\cdot \GRAD_x g=\frac12 
w\wedge\big(3(w.B)(B\wedge\GRAD_B B)-B\wedge
\GRAD_w B-\GRAD_{B\wedge w} B\big)\cdot \GRAD_w g$$
\end{itemize}
with the notation $ \GRAD_V \Phi \eqdefa V \cdot \GRAD   \Phi $.
}
\end{Thm}
The result obtained in this paper is very different because of the
incompressibility constraint, which imposes a lot of rigidity to the system. In
particular, the kernel of the penalization is much smaller and the 
limiting system has
less degrees of freedom.

\subsubsection{A remark in the inviscid case}
  \label{Euler} The weak compactness
method used here allows to study the singular limit without 
regularity with respect to
the time variable. However it uses crucially the strong compactness 
in $x$ given by the
energy estimate (\ref{eq:energie}). Implicitly we have actually considered the
penalization
$$L_\eps:u\in \H \mapsto P(u\wedge B)-\eps \LAP u \in H^{-1}(\Omega).$$

That rules out the possibility to manage an analogous study for 
inviscid rotating
fluids, the first obstacle being  to prove the existence of solutions for
$$\d_t \ue +(\ue\DOT \GRAD) \ue +\frac1\eps  \ue\wedge B +\GRAD p=0, 
\quad \DIV \ue=0$$
on a uniform time interval $[0,T]$. Indeed the operator
$\exp(tL/\eps)$ is  not even
uniformly bounded on $H^s(\Omega)$ for $s\geq \frac12 \cdotp$

\begin{Lem} {\sl
Define the linear operator $L$ by (\ref{L-def}). Then, the group
$(\exp(tL))_{t\in \R}$ generated by
$L$ is not uniformly bounded on $H^s(\Omega)$ for $s\geq \frac12 \cdotp$}
\end{Lem}

\begin{proof}
The proof of that result is simply due to the fact that by definition
of~$ \Pi$ seen above, the trace of~$ \Pi u$ is not defined on~$
\partial {\mathcal S}$ even if~$ u \in \H$. So~$ \Pi$ is not continuous
on~$H^s(\Omega)$ for $s\geq \frac12 \cdotp$ Then the formula
$$
\Pi=\lim_{T\to +\infty}{1\over T} \int_0^T e^{tL} dt
$$
implies that $e^{tL}$ cannot be uniformly bounded in $H^{s}(\Omega)$ for $s\geq
\frac12$ which proves the lemma.
\end{proof}



\section{The case of a constant vector field $B$~:\\ the 2D
Navier-Stokes limit}

\label{constant}
\setcounter{equation}{0}
In the previous section, we have obtained a constraint
equation on the limiting velocity field, which expresses that
$\ul$ belongs to the kernel of the singular perturbation
$L$. This comes from the fact that the component $\we=\Pp
\ue$ has fast oscillations with respect to time, and
consequently converges weakly to 0. In the case where
$\Omega_3=\R$, this characterizes completely the weak limit
$\ul=0$.

Then it remains to get an evolution equation for $\ul$ in the
case where $\Omega_3=\T$. A natural idea consists in
projecting the evolution equation (\ref{NS}) for $\ue$ on the
kernel of
$L$, and to study its limit as $\eps \to 0$. The difficulty
is to take limits in the nonlinear terms~: as
Corollary \ref{t-reg} provides strong compactness on the
non-oscillating component $\ub$, the problem comes actually
   to prove that the oscillating terms $\we$ do not
product any constructive interference.

In order to have a good understanding of this phenomena and
of its mathematical formulation, we propose to consider first
the case where the vector field $B$ is constant and
homogeneous. The convergence result established here is not
so precise as the ones given in~\cite{basil}-\cite{basil3},~\cite{I3} 
or~\cite{Grenier},
since it
does not describe the oscillating
component and consequently does not provide any strong
convergence. Nevertheless the proof is less technical (in particular 
it does not
require any knowledge on the spectral structure of~$ L$), which
allows to consider more general cases in the sequel.

\subsection{Projection on the kernel}
In order to obtain the evolution of the limiting velocity
field $\ul$, the idea is to use  the strong
convergence
$$
\ub=\Pi \ue \to \ul \hbox{ in }
L^2_{loc}(\R^+\times
(\OO \cup \SS)),
$$
rather than the weak convergence
$$
   \ue \rightharpoonup \ul \hbox{ in }
\wL^2_{loc}(\R^+\times
\Omega).
$$
Having this idea in view, we first determine the evolution
equation for $\ub$.
\begin{Prop}\label{NS-moyennee}
{\sl
Let $u^0$ be any vector field in $\V$.
For all $\eps>0$, denote by $\ue$ a weak solution of
(\ref{NS}), and by $\ub=\Pi \ue$ its projection on
$\Ker(L)$. Then,
\begin{equation}
\label{moy2}
\d _t \ub-\Pi ( \ue\wedge \ROT \ue)-\nu \LAP \ub=0.
\end{equation}
}
\end{Prop}

\begin{proof} Identity (\ref{moy2}) is essentially a
variant of  (\ref{moy}). Indeed, in the case of a constant
$B$, the projection $\Pi$ commutes with any partial
derivative, in particular with the Laplacian $\LAP$~:
$$\Pi (\nu \LAP \ue) =\nu \LAP \Pi \ue =\nu \LAP \ub.$$
Then the key argument is the following identity~:
\begin{equation}
\label{conv-identity}
u\wedge \ROT u =\GRAD {|u|^2\over 2} -\DIV(u\otimes u)+u\DIV u
\end{equation}
As $\ue$ is divergence-free,
$$ \ue \DIV \ue=0.$$
As $\Ker(L) $ is embedded in the space of divergence-free
vector fields,
$$\Pi \left(\GRAD {|\ue|^2\over 2}\right) =0.$$
Replacing in (\ref{moy}) leads to the expected result.
\end{proof}

\subsection{ Brief description of the oscillations}
\begin{Prop}\label{propag}
{\sl
Let $u^0$ be any vector field in $\V$.
For all $\eps>0$, denote by $\ue$ a weak solution of
(\ref{NS}), and by $\we=\Pp \ue$ its projection on
$\Ker(L)^\perp$. Then there exists $\We\in
L^2_{loc}(\R^+\times \Omega)$ such that $\we=\d _3 \We$.
Moreover,
\begin{equation}
\label{oscil-eq}
\begin{aligned}
\eps &\d _t (\d_3 \We-\GRAD \Wet)+\d_3 \We\wedge B = r_\eps,\\
\eps &\d _t (\d_2 \Weu -\d_1 \Wed) +b(\d _1 \Weu +\d _2
\Wed)=s_\eps,
\end{aligned}
\end{equation}
where $r_\eps, s_\eps$ converge to 0 in
$L^2_{loc}(\R^+,H^{-5/2}_{loc}(\OO \cup \SS))
$ as
$\eps \to 0$.
}
\end{Prop}
\begin{proof} As we have supposed in this section that $\SS=\Omega$, the
projection of $\ue$ on $\Ker(L)$ satisfies, due to 
Proposition~\ref{decomposition},
$$\Pi \ue= {1\over |\Omega_3|}\int \ue dx_3.$$
Then
$$\int \we dx_3=\int (\ue -\Pi\ue) dx_3 =0,$$
and there exists $\We$ such that
$$\we=\d_3\We.$$
Moreover, as $\Omega_3=\T$, we can always choose $\We$ so that
$$\int \We dx_3=0.$$

It remains then to determine the equations for $\We$.
Equation (\ref{NS}) implies
$$\eps \d _t \we+\we\wedge be_3 +\GRAD p=\eps \nu\Pp \LAP \ue-\eps
\Pp(\ue\DOT \GRAD \ue),$$
which can be rewritten in terms of $\We$
\begin{equation}
\label{We-eq}
\eps \d _t \d_3 \We+\d_3 \We\wedge be_3 +\GRAD p=S_\eps ,
\end{equation}
with
$$S_\eps=\eps (\nu \Pp\LAP
\ue-
\Pp(\ue\DOT
\GRAD \ue)).$$
  From the energy bound, we deduce that the right hand side in (\ref{We-eq}) is
of order
$\eps$ in  the space~$L^2_{loc}(\R^+,H^{-3/2}_{loc}(\OO \cup \SS))
$. Indeed, using the continuity properties of~$ \Pp$ as well
as~(\ref{eq:inclusion}) we get
$$
\|\Pp \LAP \ue \|_{L^2([0,T], H^{-1}(\OO \cup \SS))}\leq
\|\LAP \ue \|_{L^2([0,T],  H^{-1}(\Omega))}\leq 
\|\ue\|_{L^2([0,T],\dot H^{1}(\Omega))}
$$
and
\begin{eqnarray*}
  \|\Pp(\ue\DOT
\GRAD \ue)\|_{L^2([0,T],H^{-3/2}(K))} & \leq  &  \|\ue\DOT
\GRAD \ue\|_{L^2([0,T],H^{-3/2}(K))} \\
& \leq & C \|\ue\|_{L^\infty(\R^+,L^2(\Omega))}
\|\GRAD \ue\|_{L^2([0,T],L^2(\Omega))}
\end{eqnarray*}
for all compact subsets $K$ of $\OO \cup \SS$.
Moreover, the right-hand side in (\ref{We-eq}) belongs to the image of~$\Pp$
and can therefore be written as a partial derivative with respect to $x_3$,
$$
S_\eps = \partial_3 R_\eps \quad \mbox{with} \quad \int R_\eps \: dx_3 = 0.
$$
The
equation on the third component provides then
$$\eps \d_t \d _3 \Wet +\d_3p=\d_3 R_\eps$$
where $R_\eps$ converges to 0 in
$L^2_{loc}(\R^+,H^{-3/2}_{loc}(\OO \cup \SS))
$ as
$\eps \to 0$.
Integrating with respect to $x_3$ provides, since $\int p \:dx_3 = 0$,
$$\eps \d_t  \Wet +p= R_\eps.$$
Replacing in (\ref{We-eq}) leads to
$$
\eps \d _t \d_3 \We+\d_3 \We\wedge be_3 +\GRAD (R_\eps-\eps
\d _t
\Wet)=S_\eps$$ which is the first identity in Proposition \ref{propag}.

In order to establish the second identity, we compute the rotational of
(\ref{We-eq}) and write its last component
$$\eps \d _t \d_3 (\d_2 \Weu-\d_1 \Wed)+\d_3 (\d_1 \Weu+\d_2\Wed)
=(\d_2 S_{\eps,1}-\d_1 S_{\eps,2}). $$
Moreover, the right hand side in (\ref{We-eq}) belongs to the image of
$\Pp$ (recall that in this section~$ \Pp $ commutes with all derivatives)
and can therefore be written as a partial derivative with respect to $x_3$.
Then integrating with respect to $x_3$ leads to the expected equality, where
the right-hand side converges to 0  in
$L^2_{loc}(\R^+,H^{-5/2}_{loc}(\OO \cup \SS))
$ as~$\eps$ goes to~$0$.
\end{proof}

\subsection{Study of the coupling}
   The algebraic structure of the propagator (\ref{oscil-eq})
implies that the oscillating terms cannot interact and produce
some contribution in the limiting equation governing  $\ul$.
Indeed  the nonlinear term
$w_\eps
\wedge \ROT w_\eps$ can be rewritten as the sum of  a total
derivative with respect to $x_3$ and a total
derivative with respect to $t/\eps$, modulo a remainder
which converges formally to 0.  Then, in order to prove a
rigorous convergence result, the first step is to introduce a
regularization of the equations (\ref{oscil-eq}) and to get a
control of the source terms in some strong norm.

\begin{Lem}\label{regularisation}
{\sl
Let $u^0$ be any vector field in $\V$.
For all $\eps>0$, denote by $\ue$ a weak solution
of~(\ref{NS}) in~$ L^\infty(\R^+,\V) \cap L^2(\R^+,\dot H^1)$,
and by~$\we=\Pp \ue$ its projection onto~$\Ker(L)^\perp$.
Then, for all~$\delta>0$, there exists~$\Wd\in
L^2_{loc}(\R^+,\cap _s H^s(\Omega))$ such that
\begin{equation}
\label{regular}
   \we-\d_3 \Wd
\to 0 \hbox{  in } L^2_{loc}(\R^+\times \Omega)\hbox{
as }\delta
\to 0
\hbox{ uniformly in }\eps>0,
\end{equation}
and
\begin{equation}
\label{oscil-app}
\begin{aligned}
\eps &\d _t (\d_3 \We^\delta -\GRAD \Wet^\delta )+\d_3 \We^\delta \wedge B =
r^\delta_\eps,\\
\eps &\d _t (\d_2 \Weu^\delta  -\d_1 \Wed^\delta ) +b(\d _1 \Weu^\delta  +\d _2
\Wed^\delta )=s^\delta_\eps,
\end{aligned}
\end{equation}
where  for all~$ \delta$, $r_\eps^\delta$ and~$ s^\delta_\eps$ converge to 0 in
$L^2_{loc}(\R^+,L^2_{loc}(\Omega))
$ as
$\eps \to 0$.
}
\end{Lem}

\begin{proof}
We introduce the following regularization: let~$\kappa \in 
C^\infty_c(\R^3,\R^+)$ such that~$\kappa(x)=0$ if~$|x|\geq 1$ 
and~$\int \kappa dx =1$, we define
$$\kappa_\delta:x\mapsto {1\over \delta^3}\kappa\left(
{.\over
\delta}\right) $$
and
$$\wde=\we *\kappa_\delta,\quad\Wd=\We* \kappa_\delta,$$
so that
$\wde =\d_3 \Wd.$

   By the energy estimate, for all $T>0$,
$$
\ue \hbox{ is uniformly bounded in
}L^2([0,T],H^1(\Omega)).
$$
  By Proposition
\ref{decomposition}, we infer that
$$\we=\Pp \ue \hbox{ is uniformly bounded in
}L^2([0,T],{\mathcal K})$$
for some compact~$ {\mathcal K}$ of~$ L^2(\Omega)$.
The result (\ref{regular}) then follows from the following fact:
$$
\wde(t,x)-\we(t,x)=\int
(\we(t,x-y)-\we(t,x))\kappa_\delta(y) dy
$$
hence there is some  continuity modulus~$ \omega$ such that
$$
\forall \varepsilon, \quad |\wde(t,x)-\we(t,x)| \leq \int \omega(y) 
\kappa_\delta(y) dy,
$$
and
the result follows.

Regularizing (\ref{oscil-eq}) leads to
$$\begin{aligned}
\eps &\d _t (\d_3 \Wd-\GRAD \Wdt)+\d_3 \Wd\wedge B =
r_\eps* \kappa_\delta,\\
\eps &\d _t (\d_2 \Wdu -\d_1 \Wdd) +b(\d _1 \Wdu +\d _2
\Wdd)=s_\eps*\kappa_\delta,
\end{aligned}
$$
because $b$ is homogeneous. Then, for all $T>0$ and all
compact subsets~$K$  of~$\OO \cup \SS$ and for $\delta$ small enough,
$$
\begin{aligned}
\|r_\eps^\delta\|_{L^2([0,T],L^2(K))}&=\|r_\eps*\kappa_\delta\|_{L^2([
0,T],L^2(K))}\\
&
\leq
C\|\kappa_\delta\|_{W^{5/2,1}(\R^3)}\|r_\eps\|
_{L^2([0,T],H^{-5/2}(K'))}\\
&\leq {C\over \delta^{5/2}}\|r_\eps\|
_{L^2([0,T],H^{-5/2}(K'))},
\end{aligned}
$$
where $K'$ is a compact subset of $\OO \cup \SS$ such that
$$\{x\in \Omega\,/\, d(x,K)\leq \delta\} \subset K'.$$
And, in the same way,
$$\|r_\eps^\delta\|_{L^2([0,T],L^2(K))}
\leq {C\over \delta^{5/2}}\|r_\eps\|
_{L^2([0,T],H^{-5/2}(K'))}.$$
For a fixed $\delta$, Proposition \ref{propag} gives  the
expected convergences.
\end{proof}

Equipped with this preliminary result, we are now able to
study the coupling between the oscillating terms and to prove
the following proposition.

\begin{Prop}\label{couplage}
{\sl
Let $u^0$ be any vector field in $\V$.
For all $\eps>0$, denote by $\ue$ a weak solution of
(\ref{NS}), and by $\we=\Pp \ue$ its projection onto~$\Ker(L)^\perp$. Then,
$$P\int \we\wedge \ROT \we dx_3 \to 0 \hbox{ in }
\DD'(\R^+\times \Omega) \hbox{ as }\eps
\to 0,$$
where $P$ denotes the Leray projection.
}
\end{Prop}

\begin{proof}
We start by proving 
that
\begin{equation}
\label{approx}
P\left(\int \we\wedge \ROT \we 
dx_3 -\int \wed
\wedge
\ROT
\wed dx_3\right)\to 0 \quad \mbox{in} 
\quad \DD'(\R^+\times \Omega) \:
  \hbox{ as }\: \delta
\to 0 \: 
\hbox{ uniformly in }\:  \eps.
\end{equation}
 From the 
identity
(\ref{conv-identity})
and the relations $\DIV 
\we=\DIV\wed=0$, we deduce that
$$
\begin{aligned}
P&\left(\we\wedge 
\ROT 
\we-\wed
\wedge
\ROT
\wed\right)\\
&=P\GRAD\left({|\we|^2-|\wed|^2\over 2}\right)-P\DIV
(\we\otimes \we-\wed\otimes \wed)\\
&=-P\DIV 
((\we-\wed)\otimes \we) 
-P\DIV(\wed\otimes
(\we-\wed))
\end{aligned}
$$
By Lemma 
\ref{regularisation}, we get (\ref{approx}).

It remains  to prove 
that, for any fixed $\delta>0$,
\begin{equation}
\label{approx2}
\int 
\wed
\wedge
\ROT
\wed dx_3\to 0
  \hbox{ as }\eps
\to 
0.
\end{equation}
Define $\rho_{\eps}^\delta=\ROT 
\Wd$~:
$$\left(\begin{array}l
\ru\\
\rd 
\\
\rt
\end{array}\right)
=
\left(\begin{array}l
  \d_2 \Wdt-\d_3 
\Wdd\\
  -\d_1 \Wdt+\d_3 \Wdu\\
  \d_1 \Wdd-\d_2 
\Wdu
\end{array}\right)
$$
which is uniformly bounded with respect to 
$\eps$ in
$L^2([0,T], \cap_s H^s(\Omega))$ by 
Lemma
\ref{regularisation}. We have
$$\wed\wedge \ROT 
\wed=\left(\begin{array}l
\d_3 \Wdu\\
\d_3 \Wdd\\
\d_3 
\Wdt
\end{array}\right)
\wedge
\left(\begin{array}l
\d_3 \ru\\
\d_3 
\rd\\
\d_3 \rt
\end{array}\right).
$$
 From (\ref{oscil-app}) and the 
divergence-free relation
$\d_3
\Wdt =-\d_1
\Wdu-\d_2 \Wdd$, we deduce 
that the previous term can be
rewritten
$$\left(\begin{array}l
\eps 
b^{-1}\d _t (\d _3 \Wdd-\d_2
\Wdt)-b^{-1}r_{\eps,2}^\delta\\
\eps 
b^{-1}\d _t (-\d _3 \Wdu+\d_1
\Wdt)+b^{-1}r_{\eps,1}^\delta\\
\eps 
b^{-1}\d _t (\d _2 \Wdu-\d _1 
\Wdd)-b^{-1}s_\eps^\delta
\end{array}\right)\wedge
\left(\begin{array} 
l
\d_3 \ru\\
\d_3 \rd\\
\d_3 \rt
\end{array}\right),
$$
or 
equivalently
$$\left(\begin{array}l
-\eps b^{-1}\d _t 
\ru
-b^{-1}r_{\eps,2}^\delta\\
-\eps b^{-1}\d _t 
\rd+b^{-1}r_{\eps,1}^\delta\\
-\eps b^{-1}\d _t 
\rt-b^{-1}s_\eps^\delta
\end{array}\right)\wedge
\left(\begin{array}l 
\d_3 \ru\\
\d_3 \rd\\
\d_3 \rt
\end{array}\right).
$$
Integrating by 
parts with respect to $x_3$ leads then to
$$b^{-1} 
\left(
\begin{array}l
-\eps \d _t (\rd \d _3 \rt)+\d _3 (\rd \eps 
\d
_t \rt) +(r_{\eps,1}^\delta \d _3 \rt+s_\eps ^\delta 
\d
_3\rd)\\
-\eps \d _t (\rt \d _3 \ru)+\d _3 (\rt \eps  \d
_t \ru) 
+(r_{\eps,2}^\delta \d _3 \rt-s_\eps ^\delta \d
_3\ru)\\
-\eps \d _t 
(\ru \d _3 \rd)+\d _3 (\ru \eps  \d
_t \rd) -(r_{\eps,2}^\delta \d _3 
\rd+r_{\eps,1} ^\delta \d
_3\ru)
\end{array}\right),$$
from which we 
deduce that, for any fixed $\delta>0$,
$$
P \int \wed \wedge \ROT 
\wed dx_3 \to 0 $$
in the sense of distributions, as $\eps \to 
0$.

Combining (\ref{approx}) and (\ref{approx2}) gives the
expected 
convergence.
\end{proof}

\subsection{Passage to the limit}
In order 
to determine the limiting velocity field $\ul$, we have now to 
take
limits in (\ref{moy2}) which can be rewritten
$$\d _t \ub-\Pi ( 
\ub\wedge \ROT \ub+\we \wedge \ROT \ub +\ub\wedge \ROT
\we +\we\wedge 
\ROT \we)-\nu
\LAP
\ub=0,
$$
using the 
decomposition
$$\ue=\ub+\we,$$
where we recall that
$$\ub \to \ul 
\hbox{ strongly in }L^2_{loc}(\R^+\times \Omega),$$
$$\we 
\rightharpoonup  0 \hbox{ weakly in }L^2_{loc}(\R^+\times 
\Omega).$$

Then standard arguments using~(\ref{conv-identity}) show 
that
$$\d_t \ul-\Pi ( \ul\wedge \ROT \ul)-\nu
\LAP
\ul=\lim_{\eps\to 
0} \Pi (\we\wedge \ROT \we),
$$
and by  Proposition \ref{couplage} we 
get
$$\d_t \ul-\Pi ( \ul\wedge \ROT \ul)-\nu
\LAP
\ul=0.
$$
In the 
case where $B$ is constant, the projection $\Pi$ reduces to a 
simple
averaging with respect to $x_3$, and we recover the usual 
convergence result
towards the 2D1/2 Navier-Stokes equation.

\section{The case of a variable vector field $B$~:\\
a turbulent behaviour}
\label{variable}
\setcounter{equation}{0}
In this section we shall prove Theorem~\ref{thm:variable} stated in
the introduction, concerning the general case when the rotation
vector~$ B = b(x_h) e_3$ is inhomogeneous. We suppose that Assumptions~$ (H0)$
to~$ (H2)$ are satisfied.
If~$ \Omega_3 = \R$, then~$ \overline u = 0$ simply because it is in~$
L^2(\Omega)$ but only depends on the horizontal variables. So from now
on we can suppose that~$\Omega_3 = \T$. 

The strategy of proof is quite similar to the constant case, so we
shall often be referring to the results of the previous
section. The first remark to be made is that if~$  b$ is
constant in some positive measured region of~$ \Omega_h$, then in that
region the results of the previous section should apply and one
should recover at the limit the usual two--dimensional behaviour.

Moreover, the results of Section~\ref{sing} hold for any~$ B$, so in
particular any weak limit point~$ \ul$ of a sequence of weak
solutions~$ \ue$ to~(\ref{NS}) is in the kernel of~$ L$ according to
Corollary~\ref{contrainte}. That means in particular that the third
component does not see the difference between~$ \OO$ and~$ \SS$ since the elements of the kernel of~$ L$ have the
same third component whether~$ b$ is homogeneous or not. So in the following, we shall restrict the study
of the limit system to the horizontal components only.
 As in the previous section, the proof
of Theorem~\ref{thm:variable} consists in finding the equation
satisfied by~$\ul$ (at least its horizontal part~$ \ul_h$), by taking
the 
limit of the equation satisfied
by the horizontal part of~$\ub = \Pi \ue$. 
The first result we shall establish is that in the
general, variable~$ b$ case, there is no coupling between oscillating
vector fields yielding extra terms in the averaged equation. This will
be a generalization of Proposition~\ref{couplage} to the variable
case, and the analysis will follow closely the proof of
Proposition~\ref{couplage}. Then we shall write the averaged equation
on~$ \OO$. Finally we shall concentrate on the~$ \SS$ case
and show the limit~$ \ul_h$ satisfies a two--dimensional Navier--Stokes
equation with homogeneous boundary conditions on the boundary of~$ \SS $.

\subsection{The averaged equation}
Let~$ \ue$ be a family of weak solutions to~(\ref{NS}), and
define~$\ub = \Pi \ue $. 
Recalling that the elements of~$ \Ker
(L)$ are divergence free, we have as in the constant case
$$
\partial_t \ub - \nu  \Pi \Delta \ue  - \Pi (\ue \wedge \ROT \ue)  = 0.
$$
Of course  the projector~$ \Pi$ no longer commutes with
(horizontal) derivatives. However as~$\Pi$ belongs to  $C(w-H^s(\OO \cup \SS),
w-H^s(\OO \cup \SS))$ for~$ s \leq 0$,
we clearly have as~$ \eps$ goes to zero,
$$
  \Pi \Delta \ue   \rightharpoonup  \Pi \Delta \ul.
$$
Let us now take the
limit in the quadratic term. 

\begin{Lem}
\label{limitquad}
{\sl
Let $u^0\in L^2(\Omega)$ be a divergence-free vector field.
For all~$\eps>0$, denote by $\ue$ a weak solution of~(\ref{NS}), and
by~$ \ub = \Pi \ue$ (resp. ~$\we=\Pp \ue$) its projection onto~$\Ker (
L) $ (resp. 
$\Ker(L)^\perp$). Then the following results hold in~$ {\mathcal
  D}'(\R^+\times (\OO \cup \SS))$, as~$ \eps$ goes
to zero:
 \begin{equation}
\label{eq:firstlimit}
\Pi (\ub \cdot \nabla \ub )  \to  \Pi (\ul \cdot \nabla \ul),
\end{equation}
\begin{equation}
\label{eq:secondlimit}
\Pi (\we  \cdot \nabla \ub  + \ub  \cdot \nabla \we ) \to   0,
\end{equation}
\begin{equation}
\label{eq:thirdlimit}
\Pi (\we  \cdot \nabla \we) \to  0.
\end{equation}
}
\end{Lem}
\begin{proof}
The results~(\ref{eq:firstlimit}) and~(\ref{eq:secondlimit})  are
simply due to the compactness of~$ \ub$ as well as the fact that~$
\we$ goes weakly to zero,  results given in Corollary~\ref{t-reg}. We
also use the continuity properties of $\Pi$ stated in Proposition~\ref{decomposition}.
The more difficult result to prove is of
course~(\ref{eq:thirdlimit}). The method will follow the proof of
Proposition~\ref{couplage}, and will be achieved in two steps. First
we show that one can smoothen out the equation satisfied by $w_\eps$, and then
we perform some algebra on the bilinear term in the equation, as in the constant case. We shall therefore continuously be referring to
the methods of Section~\ref{constant}.

Let us start by proving the following result, analogous to
Proposition~\ref{propag}.
\begin{Prop}
\label{propagvariable}
{\sl
Let $u^0\in L^2(\Omega)$ be a divergence-free vector field.
For all~$\eps>0$, denote by~$\ue$ a weak solution of
(\ref{NS}), and by $\we=\Pp \ue$ its projection
onto~$\Ker(L)^\perp$. Then there exists~$ W_{\eps,3} \in
L^2_{loc}(\R^+\times \Omega)$ such that $w_{\eps,3}=\d_3  W_{\eps,3}$.
Moreover,
\begin{equation}
\label{oscil-eqvariable}
\begin{aligned}
\eps &\d _t (w_{\eps,h} - \nabla_h \WWe) + w_{\eps,h} \wedge B = \tilde
r_\eps \\
\eps &\d _t (\d_2 \weu -\d_1 \wedd) + \mbox {div}_h (b w_{\eps,h}  )=\tilde s_\eps,
\end{aligned}
\end{equation}
where~$\tilde r_\eps, \tilde s_\eps$ converge to 0 in
$L^2_{loc}(\R^+,H^{-5/2}_{loc}(\OO \cup \SS))
$ as
$\eps \to 0$. 
}
\end{Prop}

\begin{proof}
We shall omit the  proof of that result here, as it is identical to
  the proof of  Proposition~\ref{propag}: we just have to notice that
  the third component~$w_{\eps,3} $ is of vertical mean zero, so can
  as in the constant case be replaced by~$\d_3  W_{\eps,3} $. The
  other components, contrary to the constant case, cannot be
  transformed in that way, so remain as they are. The rest of the
  proof is identical to the constant case.
\end{proof}

Now as in the constant case, let us smoothen out
Equation~(\ref{oscil-eqvariable}). 
\begin{Lem}
{\sl
Let $u^0\in L^2(\Omega)$ be a divergence-free vector field.
For all~$\eps>0$, denote by~$\ue$ a weak solution of
(\ref{NS}), and by~$\we=\Pp \ue$ its projection on
$\Ker(L)^\perp$. Then, for all $\delta>0$, there
exists
$w_{\eps,h}^\delta\in L^2_{loc}(\R^+,\cap _s H^s(\Omega))$ and~$ \WWd \in
L^2_{loc}(\R^+,\cap _s H^s(\Omega))$ such that
\begin{equation}
\label{regularvariable1}
  w_{\eps,h}-w_{\eps,h}^\delta
\to 0 \hbox{  in } L^2_{loc}(\R^+\times \Omega)\hbox{
as }\delta
\to 0
\hbox{ uniformly in }\eps>0,
\end{equation}
\begin{equation}
\label{regularvariable2}
  w_{\eps,3}- \partial_3 \WWd
\to 0 \hbox{  in } L^2_{loc}(\R^+\times \Omega)\hbox{
as }\delta
\to 0
\hbox{ uniformly in }\eps>0,
\end{equation}
\begin{equation}
\label{gradient}
\delta \|\nabla_h w_{\eps,h}^\delta\|_{L^2(\Omega)} \to 0 \hbox{
as }\delta
\to 0
\hbox{ uniformly in }\eps>0
\end{equation}
and
\begin{equation}
\label{oscil-appvariable}
\begin{aligned}
\eps &\d _t (w_{\eps,h}^\delta - \nabla_h \WWd) +  w_{\eps,h}^\delta \wedge B = \tilde
r^\delta_\eps \\
\eps &\d _t (\d_2 \wdu -\d_1 \wdd) + \d_1 (\wdu b) + \d_2 (\wdd b)
=\tilde s^\delta_\eps,
\end{aligned}
\end{equation}
where  $\tilde r_\eps^\delta $ and $\tilde s^\delta_\eps$ converge to 0 in
$L^2_{loc}(\R^+,L^2_{loc}(\OO \cup \SS))
$ as
$\eps$ and~$ \delta$ go to~$0$. More precisely  for any subset~$ K$
of~$ \OO \cup \SS$ there is a constant~$ C$ (independent of~$ \e$ and~$
\delta$) and a constant~$ c_\delta$ depending only on~$ \delta$ such that
\begin{equation}
\label{estimateremainder}
\|\tilde r_\eps^\delta + \tilde s^\delta_\eps\|_{L^2 (K)} \leq 
c_\delta o_\eps(1) + C \delta.
\end{equation}
}
\end{Lem}
\begin{proof}
We shall not write all the details of the proof here, since it is
very similar to the constant case (Lemma~\ref{regularisation}); let us
simply point  out where the fact
that~$ b$ is not constant appears --- note that
in~(\ref{estimateremainder}), the part~$c_\delta o_\eps(1) $ is
precisely due to the terms of the constant case, and we will see here
that the fact that~$ b$ is no longer constant yields terms which are
estimated by~$ C \delta $. In the approximation of the
equation, the only difference with the
constant case is that of course~$ (w_{\eps} \wedge B) * \kappa_\delta$
is not equal to~$  (w_{\eps} * \kappa_\delta) \wedge B$. Moreover of
course~(\ref{gradient}) is obvious in the constant case, since $\Pi_\perp$ commutes with partial derivatives. So we need to
deal with those two problems due to the fact that~$ b$ is not constant.

First of all, the
difference between~$ (w_{\eps} \wedge B) * \kappa_\delta$ and~$  (w_{\eps} * \kappa_\delta) \wedge B$ is small when~$\delta$ goes to zero, due to
the following computation: we have
$$
 (w_{\eps} \wedge B) * \kappa_\delta (t,x) -  (w_{\eps} * \kappa_\delta)
 \wedge B(t,x) = \int  w_{\eps} (t,x-y)  \kappa_\delta (y)   \wedge ( B
 (x-y) - B(x)) \: dy
$$
hence
$$
(w_{\eps} \wedge B) * \kappa_\delta (t,x) -  (w_{\eps} * \kappa_\delta)
 \wedge B(t,x) = \int_0^1 \int y \cdot \nabla B (x - \sigma y)\wedge   w_{t,\eps}  (x-y)  \kappa_\delta (y) \: dy d\sigma.
$$
To conclude we need to take the~$ L^2$ norm in~$x$ of that quantity,
and Young's inequality yields
$$
\|
(w_{\eps} \wedge B) * \kappa_\delta (t,\cdot ) -  (w_{\eps} * \kappa_\delta)
 \wedge B(t,\cdot) \|_{L^2(\Omega) } \leq C_B \||\cdot|  \kappa_\delta
 \|_{L^1(\Omega) } \|  w_{\eps}(t)\|_{L^2(\Omega) }  \leq C \delta \|u_0\|_{L^2}
$$
uniformly in time. The result follows for the first equation
in~(\ref{oscil-appvariable}). The second one is of the same type,
since
$$
\d_1 (w_{\eps,1} b )  + \d_2 (w_{\eps,2} b ) = b \mbox{div}_h
w_{\eps,h}  + w_{\eps,h} \cdot \nabla_h b.
$$
The term~$(w_{\eps,h} \cdot \nabla_h b) *   \kappa_\delta $ is
approximated by~$( w_{\eps,h} * \kappa_\delta)  \cdot \nabla_h b$
exactly as above; to replace the term~$( b \mbox{div}_h
w_{\eps,h}) *  \kappa_\delta  $ by~$b ( \mbox{div}_h
w_{\eps,h} *  \kappa_\delta   )$ we write the same type of
computation, with
$$
( b \mbox{div}_h
w_{\eps,h}) *  \kappa_\delta (t,x)  - b ( \mbox{div}_h
w_{\eps,h} *  \kappa_\delta   )(t,x)  = - \int   \partial_3 w_{\eps,3} (t,x-y)  \kappa_\delta (y)   ( b
 (x-y) - b(x)) \: dy
$$
hence, since $\Pi_\perp$ commutes with $\partial_3$,
$$
\|( b \mbox{div}_h
w_{\eps,h}) *  \kappa_\delta  - b ( \mbox{div}_h
w_{\eps,h} *  \kappa_\delta   ) \|_{L^2(\R^+\times \Omega) } \leq C_B
\|\partial_3 w_{\eps,3}\|_{L^2(\R^+\times \Omega)} \delta \leq  C_B \delta \|u_0\|_{L^2}.
$$
That ends the proof of~(\ref{oscil-appvariable}). 

Now to end the proof of the proposition, we still need to
check~(\ref{gradient}). The idea is to use the following estimate, due
to the fact that~$ w_\e$ is bounded in~$ L^2_t({\mathcal K})$ for some
compact subspace of~$ L^2(\Omega)$: there is a
continuity modulus~$ \omega$ such that
\begin{equation}
\label{eq:continuitymod}
\forall y \in \Omega, \quad \|w_\e (t,\cdot+y) - w_\e (t,\cdot )\|_{L^2(\Omega)}
\leq \omega (y), \quad \mbox{uniformly in }\: t \: \mbox{and } \: \e. 
\end{equation}
Now since~$w_{\eps,h}^\delta  = \kappa_\delta * w_{\eps,h}$, we have
$$
\nabla_h w_{\eps,h}^\delta  (x) = \int_\Omega 
\frac{1}{\delta^4} (\nabla_h \kappa) (\frac{y}{\delta}) w_{\eps,h}(x-y) \: dy.
$$
Since~$ \int_\Omega \nabla_h \kappa (y)\: dy = 0 $, it follows that
$$
\nabla_h w_{\eps,h}^\delta  (x) = \int_\Omega 
\frac{1}{\delta^4} (\nabla_h \kappa) (\frac{y}{\delta}) 
(w_{\eps,h}(x-y) - w_{\eps,h}(x) )\: dy.
$$
Then by~(\ref{eq:continuitymod}) we find that
$$
\|\nabla_h w_{\eps,h}^\delta \|_{L^2(\Omega)} \leq C \frac{
  \eta (\delta)}{\delta} \|\nabla_h \kappa \|_{L^1(\Omega)}
$$
where~$  \eta (\delta) $ goes to zero as~$ \delta$ goes to
zero. The result is proved.

\end{proof}
Now we are ready to prove the following result.
\begin{Prop}
\label{couplagevariable}
{\sl
Let $u^0\in L^2(\Omega)$ be a divergence-free vector field.
For all~$\eps>0$, denote by~$\ue$ a weak solution of~(\ref{NS}), and
by $\we= \Pp \ue$ its projection onto~$\Ker(L)^\perp$. Then
$$
\Pi \left( \we\wedge \ROT \we \right) \to 0 \hbox{ in }
\DD'(\R^+\times (\OO \cup \SS)) \hbox{ as }\eps
\to 0.$$
}
\end{Prop}
\begin{proof}
Since the result we are looking for is a weak limit, we can
restrict our attention to the set~$ \OO$ where~$ \nabla b$ does not
vanish, as in~$\SS $ the result is due to
Proposition~\ref{couplage}. Moreover, to prove the result one can
restrict our attention to~$
\int \we\wedge \ROT \we \: dx_3
$ and it is enough to prove that it is proportionnal to~$ \nabla_h b$
(up to a small remainder term):
taking the scalar product with a function in~$\Ker(L) $ will then
yield automatically zero by definition of~~$\Ker(L)$.

In order to simplify the computations, we shall directly prove the
result replacing~$ w_{\eps,h}$ by~$w_{\eps,h}^\delta$ and~$ \wet$
by~$  \partial_3 \WWd$. The difference in the two computations is
indeed small when~$ \delta$ is small, uniformly in~$ \eps$, exactly as
in the proof of~(\ref{approx}) in the constant case. So
writing~$\partial_3 \WWd =  w_{\eps,3}^\delta $, and dropping the
index~$  \delta$ to
simplify, we can perform the following algebraic computations, which
will prove the result.

Let us start by recalling that, due to
Proposition~\ref{propagvariable}, we have
\begin{equation}
\label{eq:wedweu}
w_{\eps,h} =- \frac{1}{b} ( \eps \partial_t \rho_{\eps,h} + r_\e),
\end{equation}
where similarly to Section~\ref{constant} we have noted
$$
\rho_{\eps,h} \eqdefa \nabla_h^\perp \WWe- w_{\eps,h}^\perp .
$$
 We shall also define
$$
\rho_{\eps,3}  \eqdefa \d_1 \wedd - \d_2 \weu.
$$
In~(\ref{eq:wedweu}) and in the following, the function~$ r_\e$ denotes a remainder
term, arbitrarily small in the space~$ L^2_{loc} (\R^+ \times (\OO \cup \SS))$.
It follows that
$$
\we\wedge \ROT \we = 
\left(\begin{array}c \displaystyle 
-\frac{1}{b} (\eps\d_t  \rho_{\eps,1}+ r_\e)\\
-\frac{1}{b} (\eps\d_t  \rho_{\eps,2}+ r_\e)\\
\d_3 \WWe
\end{array}\right)
\wedge
\left(\begin{array}c
\d_3 \rho_{\eps,1} \\
\d_3 \rho_{\eps,2} \\
 \rho_{\eps,3}  
\end{array}\right),
$$
which implies that
$$
\we\wedge \ROT \we = 
\left(\begin{array}l
-\frac{1}{b} (\eps\d_t \rho_{\eps,2}+ r_\e) \rho_{\eps,3}   - \d_3 \WWe \d_3
\rho_{\eps,2} \\
\frac{1}{b} (\eps\d_t \rho_{\eps,1}+ r_\e)\rho_{\eps,3}    + \d_3 \WWe \d_3
\rho_{\eps,1}\\
\frac{1}{b}  (\eps\d_t \rho_{\eps,2} + r_\e)\d_3 \rho_{\eps,1} - \frac{1}{b}
(\eps\d_t \rho_{\eps,1} + r_\e) \d_3 \rho_{\eps,2}
\end{array}\right).
$$
Since the vertical component can be treated exactly as in the constant case, we shall now restrict our attention to the first
two components. So calling~$\alpha_h$ the horizontal components of that
vector field, we have after an integration by parts and using the
divergence free  condition~$ \partial_{33}  \WWe = - \mbox{div}_h w_{\eps,h} $:
$$
\int \alpha_h \: dx_3 = - \frac{1}{b} \int (\eps\d_t \rho_{\eps,h}^\perp+ r_\e)
\rho_{\eps,3}    \: dx_3  - \int \rho_{\eps,h}^\perp  \mbox{div}_h
w_{\eps,h}   \: dx_3  .
$$
Now we recall (calling  once again generically~$ r_\e$ the small
remainder  terms) that
$$
- \eps\d_t \rho_{\eps,3}   +  \mbox{div}_h (b w_{\eps,h})  = r_\e
$$
so
$$
  \mbox{div}_h w_{\eps,h}  =  \frac{1}{b} (
\eps \d_t \rho_{\eps,3}   -w_{\eps,h} \cdot \nabla_h
  b + r_\e)  .
$$
It follows that
\begin{eqnarray*}
\int \alpha_h \: dx_3 & = &- \frac{1}{b} \int \eps\d_t \rho_{\eps,h}^\perp
\rho_{\eps,3}    \: dx_3 -\frac{1}{b} \int \rho_{\eps,h}^\perp\eps \d_t
\rho_{\eps,3}    \: dx_3 +  \frac{1}{b} \int \rho_{\eps,h}^\perp w_{\eps,h}
\cdot \nabla_h b \: dx_3  + \int \widetilde r_\e  \: dx_3\\
 & = & -\frac{1}{b} \int \eps\d_t (\rho_{\eps,h}^\perp\rho_{\eps,3} )   \: dx_3   +  \frac{1}{b} \int \rho_{\eps,h}^\perp w_{\eps,h}
\cdot \nabla_h b \: dx_3  + \int \widetilde r_\e  \: dx_3,
\end{eqnarray*}
where now~$ \widetilde r_\e$ denotes generically  the product of~$
r_\e $ by a component of~$ \rho_\e$. But~$ \rho_{\e,3} $ is a
combination of derivatives of~$ w_\e$ whereas~$ \rho_{\e,h}$ is a
combination of components of~$ w_\e$. A product of the type~$
\rho_{\e,h}  r_\e$ clearly goes to zero in~$  {\mathcal D}'(\R^+
\times (\OO\cup \SS))$. For the term~$ r_\e \rho_{\e,3}$, one uses
result~(\ref{gradient}) and~(\ref{estimateremainder}) to infer that for any
subset~$ K$ of~$ \OO \cup \SS$,
\begin{eqnarray*}
\| r_\e \rho_{\e,3} \|_{L^1(K)} & \leq & \| r_\e \|_{L^2(K)}
\|  \rho_{\e,3}  \|_{L^2(K)} \\
 & \leq & (c_\delta o_\eps(1) + C\delta) \|\nabla_h
w_{\e,h}^\delta\|_{L^2(K)} \to 0,
\end{eqnarray*}
as $\e$ followed by $\delta$ go to zero.
So from now on~$ \widetilde r_\e$  will denote
generically a term going to zero in~$  {\mathcal D}'(\R^+
\times (\OO \cup \SS))$. 

Recalling~(\ref{eq:wedweu}) we get
$$
\int \alpha_h \: dx_3 = - \frac{1}{b} \int \eps\d_t
(\rho_{\eps,h}^\perp\rho_{\eps,3}) \: dx_3   -  \frac{1}{b^2}  \int
\rho_{\eps,h}^\perp  \eps\d_t  \rho_{\eps,h} \cdot \nabla_h b \: dx_3
+ \int \widetilde r_\e  \: dx_3 .
$$
In particular we have
$$
\int \alpha_1 \: dx_3 = - \frac{1}{b} \int \eps\d_t( \rho_{\eps,2}
\rho_{\eps,3}   ) \: dx_3 -  \frac{1}{b^2} \int \rho_{\eps,2}
\eps\d_t \rho_{\eps,1}  \d_1 b   \: dx_3-   \frac{1}{2 b^2} \int  \eps\d_t
(\rho_{\eps,2}^2 \d_2 b)  \: dx_3 + \int \widetilde r_\e  \: dx_3 .
$$
Similarly
$$
\int \alpha_2 \: dx_3 =   \frac{1}{b} \int \eps\d_t( \rho_{\eps,1}
\rho_{\eps,3}   ) \: dx_3 +  \frac{1}{b^2}   \int \rho_{\eps,1}
\eps\d_t \rho_{\eps,2} \d_2 b  \: dx_3 + \frac{1}{2 b^2}\int  \eps\d_t (\rho_{\eps,1}^2 \d_1 b)  \: dx_3 + \int \widetilde r_\e  \: dx_3 .
$$
But one can also write
$$
\int \alpha_2 \: dx_3 =   \frac{1}{b}   \int\eps\d_t ( \rho_{\eps,1}
\rho_{\eps,3} +   \frac{1}{2 b} \rho_{\eps,1}^2 \d_1 b +
\frac{1}{b}\rho_{\eps,1}\rho_{ \eps,2} \d_2 b
)  \: dx_3 -   \frac{1}{b^2} \int  \rho_{\eps,2}
\eps\d_t\rho_{\eps,1}\d_2 b \: dx_3  + \int \widetilde r_\e  \: dx_3 .
$$
It follows that up to full derivatives of the type~$\displaystyle \e \partial_t $,
$\displaystyle \int \alpha_h \: dx_3 $ is equal to
\begin{equation}
\label{eq:alphah}
\left(- \frac{1}{b^2} \int \rho_{\eps,2} \eps\d_t\rho_{\eps,1}  \:
dx_3  \right) \nabla_h b + \int \widetilde r_\e  \: dx_3 .
\end{equation}
Now the proof is almost finished: we recall that we want to take the
projector of the term~$\displaystyle \int  \we\wedge \ROT \we \: dx_3
$ onto the kernel  of~$
L$ restricted to the set~$ \OO$. Recalling that~$\Ker L$ is made of
vector fields of the type~$ \nabla_h^\perp \varphi$ with~$ \displaystyle  \nabla_h b \cdot
\nabla_h^\perp \varphi = 0$, we have obviously
$$
\Pi \left(- \frac{1}{b^2} \int \rho_{\eps,2} \eps\d_t\rho_{\eps,1}  \:
dx_3  \right) \nabla_h b = 0.
$$
The result is therefore proved for the horizontal component
of~$\displaystyle  \int
\we\wedge \ROT \we \: dx_3 $. The third component is identical to the
constant case, so the proposition is proved.
\end{proof}
This ends the proof of~(\ref{eq:thirdlimit}), hence of Lemma~\ref{limitquad}.
\end{proof}

In the following we shall denote by~$ (\LL) $ the limiting system:
$$
\displaystyle
(\LL)  \quad \quad \left\{\begin{array}{c}
\partial_t \overline u - \nu \Pi \Delta_h \overline u + \Pi (\overline u
\cdot \nabla_h \overline u ) = 0 \quad \mbox{in} \quad {\mathcal
  D}'(\R^+ \times (\OO \cup \SS)), \quad \Pp  \overline u = 0\\
\overline u_{|t = 0} = \Pi u^0.
\end{array}
 \right.
$$
The existence of solutions to this system is easy to prove, as it is of the
same form as a 2D Navier--Stokes equation: the only point which might
be a problem is that~$ \Pi$ does not in general commute with the
Friedrich frequency truncation~$ J_n$, recalled in Section~\ref{sing}: however
approximating the system by for instance
$$
\partial_t J_n \overline u_n - \nu  J_n \Pi  J_n \Delta_h J_n \overline
u_n  + J_n \Pi ( J_n\overline u_n
\cdot \nabla_h  J_n\overline u_n ) = 0, \quad   J_n \Pp  J_n  \overline u_n = 0
$$
will do the job. In any case in the next two
sections we shall give precise formulations for the solution
of~$(\LL)$: in~$ \SS$ $\overline u_h$ is the unique solution of~$ (NS2D)$, and
outside~$ \SS$ it is the solution of a heat equation.  The uniqueness
of~$ \overline u$ implies in particular that the convergence holds for
the whole sequence~$ u_\e$ and not only for a subsequence.

As noted earlier
in this section, the third component of~$\Pi u$ for any vector field~$
u$ is the same whether~$ b$ is constant or not and is simply the
vertical average of~$ u_3$. It follows that the  third component of this equation is simply
$$
\displaystyle
\partial_t  \overline u_3 - \nu \Delta_h  \overline u_3 +  \overline
u_h \cdot \nabla_h   \overline u_3 = 0 \quad \mbox{in} \quad \R^+ \times (\OO \cup \SS),
\quad  \overline u_{3|t = 0} = \int u_{3|t = 0} \: dx_3.
$$
Now all the work consists in determining~$\overline
u_h $. We shall
 consider separately the vector field on~$ \OO$  and on~$ \SS$, which is the
object of the two following sections; so in those sections, our attention
will be restricted to the
horizontal component~$ \overline u_h$.
 
\subsection{The averaged equation on $\OO$}
 \label{oo}
  We shall prove the following result, which yields the part of
  Theorem~\ref{thm:variable} which lies in~$ \OO$.
\begin{Prop}
\label{limoo}
{\sl
Let $\overline u$ be a vector field satisfying~$ (\LL)$ with
$$
\overline u \in L^\infty(\R^+,\V)\cap 
  L^2(\R^+,\dot H^1) \quad \mbox{and} \quad \overline u \in
  \mbox{Ker} (L). 
$$
Then the vector field~$ \overline u_h$ satisfies the following heat equation:
\begin{eqnarray*}
\partial_t \overline u_h - \nu \Pi \Delta_h \overline u_h =0 \quad \mbox{in} \quad
\R^+ \times \OO\\
\overline u_{h|t = 0} = \Pi u_{h,0}.
\end{eqnarray*}
}
\end{Prop}

\begin{proof}

The function~$ \overline u_h$ satisfies
$$
\partial_t  \overline u_h - \nu  \Pi \Delta_h \overline u_h + \Pi (
\overline u_h\cdot \nabla_h \overline u_h) = 0,
$$
and since~$ \overline u$ is in~$ \Ker (L)$, we have
$$
- \left( \Pi \Delta_h \overline u_h |  \overline u_h\right)_{L^2(\OO)}
= - \left(\Delta_h \overline u_h |  \overline u_h\right)_{L^2(\OO)} = 
\|\nabla_h \overline u_h\|_{L^2(\OO)}^2.
$$
We note that~$ \overline u_h$ is equal to zero on the boundary of~$
\OO$, since it is a multiple of~$ \nabla_h^\perp b$.

So to prove the proposition, the only point we need to check is that
\begin{equation}
\label{zeroL2}
\quad \forall \Phi \in \Ker (L) \cap \dot H^s(\OO), \: s > 1, 
\quad \left( \overline u_h \cdot \nabla_h
\overline u_h | \Phi_h\right)_{L^2(\OO)} = 0.
\end{equation}
By definition of~$ \Ker (L)$, we have
$$
\Phi_h \cdot \nabla_h b = \overline u_h  \cdot \nabla_h b = 0,
$$
from which we infer that
\begin{equation}
\label{parallel}
\Phi_h \wedge  \overline u_h  = 0.
\end{equation}
Now as in~(\ref{conv-identity}) we can write
$$
\Pi ( \overline u_h \cdot \nabla_h
\overline u_h) = -\Pi ( \overline u_h \wedge \mbox{rot}  \overline u_h )
$$
hence
$$
\left( \overline u_h \cdot \nabla_h
\overline u_h | \Phi_h\right)_{L^2(\OO)} = \left( \overline u_h 
\wedge  \Phi_h | \mbox{rot}  \overline u_h \right)_{L^2(\OO)} .
$$
Finally Identity~(\ref{parallel}) yields~(\ref{zeroL2}), and
Proposition~\ref{limoo} is proved.
\end{proof}

\subsection{The 2D Navier-Stokes limit on $\SS$}
In this section we shall analyse the equation satisfied by the limit
on~$\SS$, that is to say in the regions
where~$ b$ is a constant.

\begin{Prop}
\label{limitomminusoo}
{\sl
Let~$\overline u$ be a vector field satisfying~$ (\LL)$ with
$$
\overline u \in L^\infty(\R^+,\V)\cap 
 L^2(\R^+,\dot H^1)\quad \mbox{and} \quad \overline u \in
  \mbox{Ker} (L).  
$$
Then  the vector field~$ \overline u_h$ satisfies the two--dimensional Navier--Stokes
equations 
in~$\SS$, with homogeneous Dirichlet boundary conditions:
$$
\begin{array}{c}
\partial_t  \overline u_h  - \nu \Delta_h  \overline u_h + \overline u_h \cdot \nabla_h  \overline u_h = -\nabla_h p 
\quad \mbox{in} \quad
\R^+ \times \SS \\
 \overline u_{h|t= 0} =  \int u_{h,0} \: dx_3 \\
 \overline u_{h|  \partial\SS
 } = 0.
\end{array}
$$
}
\end{Prop}
\begin{proof}
First let us recall why the
equation on~$ \overline u_h $ is the two--dimensional Navier--Stokes
equation: we simply consider the weak formulation of the original
rotating fluid equations and take its limit, by integrating against a
test function~$ \Phi$, divergence free and compactly supported in~$
\R^+ \times \SS$. 
The weak formulation is as follows:
\begin{equation}
\label{eq:equepsPhi}
\int_{\R^+} \int_{\Omega} \left(-\ue \cdot \partial_t \Phi + \nu \nabla  \ue \cdot\nabla
  \Phi - \ue \otimes \ue \cdot \nabla \Phi \right)\: dx dt =
\int_{\Omega} u^0 \Phi_{|t= 0} \: dx .
\end{equation}
Then
 taking the limit as~$ \varepsilon$ goes to zero in~(\ref{eq:equepsPhi}) yields, due to Lemma~\ref{limitquad},
$$
\int \left(-\overline u \cdot \partial_t \Phi + \nu \nabla  \overline u  \cdot\nabla
  \Phi -\overline u  \otimes \overline u   \cdot \nabla \Phi \right)\: dx = \int_{\Omega} u^0 \Phi_{|t= 0} \: dx 
$$
where we have noticed that on~$\SS $, we have
$$
\Pi (\overline u \cdot \nabla \overline u  ) = \DIV (\overline u
\otimes \overline u).
$$
Now recalling that~$ \overline u$ only depends on the horizontal
variable, we deduce the expected equation on~$\overline u$, up to the
boundary terms. 
To get the boundary terms, we simply recall that the limit~$ \overline
u$ is in~$ L^2_{loc}(\R^+, H^{1}(\Omega_h))$ hence cannot have a jump on
the boundary of~$ \SS$.  Then we notice that~$ \partial \SS \subset
\partial \OO$, simply because if~$ x \in \partial\SS$, then~$ \nabla
b(x) = 0$. So  the result follows directly: the boundary condition
is a homogeneous Dirichlet boundary condition.

 Theorem~\ref{thm:variable} is proved.

\end{proof}


\end{document}